\documentclass[11pt,a4paper,reqno]{amsart}

\usepackage{times} 
\usepackage{amsmath} 
\usepackage{amssymb}  
\usepackage{amsthm}
\usepackage{latexsym}
\usepackage{amsfonts,bbm}
\usepackage{xcolor}
\usepackage{mathtools}
\usepackage{enumerate}
\usepackage{cite}
\usepackage{tikz}
\usepackage{graphicx}
\usepackage{MnSymbol} 

\newtheorem{thm}{Theorem}[section]

\newtheorem{lem}[thm]{Lemma}
\newtheorem{prop}[thm]{Proposition}
\theoremstyle{definition}
\newtheorem{defn}[thm]{Definition}
\newtheorem{rem}[thm]{Remark}

\numberwithin{equation}{section}
\allowdisplaybreaks

\renewcommand{\div}{\operatorname{div}}
\DeclareMathOperator{\tr}{tr}

\newcommand{\R}{\ensuremath{\mathbb R}}    

\newcommand{\dom}{\operatorname{dom}}
\newcommand{\setdef}[2]{\left\{ #1 \left\vert\vphantom{#1} #2 \right.\right\}}
\newcommand{\spvek}[2]{\left(\begin{smallmatrix}#1\\#2\end{smallmatrix}\right)}
\newcommand{\spvekk}[3]{\left(\begin{smallmatrix}#1\\#2\\#3\end{smallmatrix}\right)}
\newcommand{\sbvek}[2]{\left[\begin{smallmatrix}#1\\#2\end{smallmatrix}\right]}
\newcommand{\sbmat}[4]{\left[\begin{smallmatrix}#1 & #2\\#3 & #4\end{smallmatrix}\right]}

\newcommand{\Id}{I}

         \newcommand{\frakA}{\mathfrak A}

\newcommand{\calH}{\mathcal H}

\newcommand{\calR}{\mathcal R}

\newcommand{\calU}{\mathcal U}

\newcommand{\calX}{\mathcal X}         
\newcommand{\calY}{\mathcal Y}         
         
\newcommand{\<}{\langle}
\renewcommand{\>}{\rangle}

\title{Port-Hamiltonian formulation of Oseen flows}

\author[T.\ Reis]{Timo Reis}
\address{{\bf T.~Reis:} Institute for Mathematics, Faculty of Mathematics and Natural Sciences, Technische Universit\"at Ilmenau, Ilmenau, Germany}
\email{timo.reis@tu-ilmenau.de}

\author[M.\ Schaller]{Manuel Schaller}
\address{{\bf M.~Schaller:} Optimization-based Control Group, Institute for Mathematics, Technische Universit\"at Ilmenau, Ilmenau, Germany}
\email{manuel.schaller@tu-ilmenau.de}

\begin{document}
\begin{abstract}
We present Oseen equations on Lipschitz domains in a port-Hamiltonian context.
Such equations arise, for instance, by linearization of the Navier-Stokes equations.
In our setup, the external port consists of the boundary traces of velocity and the normal component of the stress tensor, and boundary control is imposed by velocity and normal stress tensor prescription at disjoint parts of the boundary. We employ the recently developed theory of port-Hamiltonian system nodes for our formulation. An illustration is provided by means of flow through a cylinder.

\end{abstract}

\maketitle

\section{Introduction}
The aim of this work is the representation of Oseen's equations in the port-Hamiltonian context. Hereby we use recently developed theory of port-Hamiltonian system nodes from \cite{PhilReis23}, which allows for an operator-theoretic description of linear infinite-dimensional port-Hamil\-tonian systems. The power of this approach is particularly due to the possibility of incorporating boundary control and observation, and it has been shown in \cite{PhilReis23} that - under suitable choice of boundary input-output configurations - Maxwell's equations, advection-diffusion equations and linear hyperbolic systems in one spatial variable fit into the framework of port-Hamiltonian system nodes.\\
Our setup consists of Oseen's equations with boundary control and observation which is formed by velocity and stress tensor in outward normal direction. Oseen's equations 
 arise in the context of linearizations around steady states of the Navier-Stokes equations, where the latter describe the dynamics of viscous Newtonian fluids. The importance of flow problems for practice and the enourmeous mathematical challenges in their treatment has led to a huge variety of textbooks and publications on their  analysis~\cite{MariTema1998,GigaNovo18,Tema95,Tema01,Raym07,Hieb20,Mon06}, numerical algorithms and finite-elements \cite{Braa07,GiraRavi12,QuarVall08} and (closed-loop) control \cite{Raym06,Fer20,Heil2022}, without claiming to deliver a complete account.
 Further, flow problems have also been considered in the port Hamiltonian context: In \cite{MoGo20}, compressible flows with boundary control are considered, where the geometric approach to port-Hamiltonian systems has been applied. In this context, incompressible flows are treated in \cite{HaMa21}, whereas Euler and shallow water flows are regarded in \cite{Swa99}. Further, an analytic approach to port-Hamiltonian formulation of flows in a one-dimensional spatial domain can be found in \cite{AltSch17}. Uncontrolled Oseen equations are considered from a Hamiltonian perspective in \cite{AcArMe22} for the purpose of a stability analysis.

This article is organized as follows: The remaining part of this introductory section is devoted to the physical derivation of Oseen's equations. Though this is somehow classical general knowledge, we give a~concise presentation in order to make this article accessible to a~community which is rather dealing with port-Hamiltonian systems than with flow problems. To motivate our configuration of boundary control and observation, we consider the flow through a~cylindric domain 
in Section~\ref{eq:cyl}. This motivates our imposed boundary conditions, which comprise inflow, normal stress prescription as well as so-called {\em do nothing} and {\em no slip conditions}. This part also contains a~formal derivation of 
the energy balance for boundary controlled Oseen's equations. Such an energy balance is typically provided by port-Hamiltonian systems. The analytical part of this article starts with Section~\ref{sec:pHnodes}, where we briefly review the theory of infinite-dimensional port-Hamiltonian systems from \cite{PhilReis23}. Thereafter, we present the main part of this article in Section~\ref{sec:oseen_node}. We consider Oseen's equations on a~bounded Lipschitz domain $\Omega\subset\R^d$, $d\ge2$. The boundary is split into two subdomains which are assumed to be separated by a~Lipschitz manifold of dimension $n-2$. These subdomains are provided with time-dependent functions with values in suitable fractional Sobolev spaces. It is in particular shown that this system fits into the framework of port-Hamiltonian system nodes.

Before starting with the derivation of the model, we declare the notation for differential operators and matrix inner products which are needed for models of flow problems. 

For a Lipschitz domain $\Omega\subset\R^d$, $d\ge2$, with boundary $\Gamma \subset \R^d$, we denote the outward unit normal $\mathrm{n} : \Gamma \to \R^d$. For the spatial variable, we use the letter $\xi=(\xi_1,\ldots,\xi_d)^\top\in\R^d$, and $t\in\R$ stands for time. By $\dot{w}$ we denote the (weak) temporal derivative of $w:I\to \mathbb{R}^{d}$, where $I\subset\R$ is an interval, and by $\nabla w = \left(\begin{smallmatrix}
\frac{\partial w}{\partial \xi_1},\ldots,\frac{\partial w}{\partial \xi_d}
\end{smallmatrix}\right)^\top$ the (weak) spatial gradient of $w:\Omega \to \mathbb{R}$. For vector-valued functions $v:\Omega \to \mathbb{R}^d$, the $\mathbb{R}^{d\times d}$-valued function $\nabla v = \left(\begin{smallmatrix}
\nabla v_1,
\ldots,
\nabla v_d
\end{smallmatrix}\right)^\top$ is the component-wise (weak) gradient, and the $\mathbb{R}^{d}$-valued function $\Delta v = \left(\begin{smallmatrix}
\Delta v_1,
\ldots,
\Delta v_d
\end{smallmatrix}\right)^\top$ is the component-wise (weak) Laplace operator. Moreover, the scalar-valued function $\div v = \sum_{i=1}^d \frac{\partial v_i}{\partial \xi_i}$ is the (weak) divergence. For a~matrix-valued functions $z=(z_{1},\ldots,z_{d})^\top:\Omega\to\R^{d\times d}$ with functions $z_{i}^\top$, $i=1,\ldots,d$, describing the rows, we denote by $\div z(\xi) = \left(\begin{smallmatrix}
\div z_{1}(\xi),
\ldots,
\div z_{d}(\xi)
\end{smallmatrix}\right)^\top$ the row-wise divergence. For vector-valued functions $b,v:\Omega\to \R^d$ we adopt the usual notation
$(b \cdot \nabla) v :=\left(\begin{smallmatrix}
b^\top \nabla v_1,\ldots, b^\top \nabla v_d
\end{smallmatrix}\right)^\top$, where $v_i:\Omega\to \R$, $i=1,\ldots,d$ are the components of $v$.
Further, for matrices $A,B\in \mathbb{R}^{d\times d}$, we will use the Hilbert-Schmidt scalar product denoted by $A:B = \tr (A^\top B)$, where $\tr$ stands for the trace of a~square matrix.

To introduce the model, let $[0,T]$ be the time interval of interest, 
let $\rho >0$ be the mass density of the fluid, and let $b : \R^3 \to\R^3$ be the flow-independent stationary convection velocity field along which the flow is linearized. The function $v:[0,T]\times \Omega\to\R^d$ stands for the velocity of the particles, $p:[0,T]\times \Omega\to\R^d$ is the infinitesimal momentum, and $P:[0,T]\times \Omega\to\R$ denotes the pressure. Velocity and infinitesimal momentum are related by the density, i.e, $p=\rho v$.

The force balance for the particles reads
\begin{subequations}\label{eq:oseen}
\begin{equation}\dot{p}(t,\xi) = \div \sigma(v(t,\xi),P(t,\xi))+ \rho\,(b(t,\xi) \cdot \nabla)\,v(t,\xi)  ,\qquad t\in[0,T],\,\xi\in \Omega, \label{eq:oseena}\end{equation}
where $\sigma(v(t,\xi),P(t,\xi))\in \R^{d \times d}$ denotes the stress tensor at $(t,\xi)$. Incompressibility of the fluid means that the velocity field is divergence-free, that is,
\begin{equation}0= \div v(t,\xi),\qquad t\in[0,T],\,\xi\in \Omega. \label{eq:oseenb}\end{equation}
Further, since $b$ is flow-independent stationary convection velocity field, we assume that its divergence as well as its normal trace vanish, that is, $\div b(\xi) = 0$ for all $\xi\in \Omega$ and $\mathrm{n}(\xi)^\top b(\xi) = 0$ for all $\xi\in \Gamma$.
The boundary of $\Omega$ is partitioned into three distinct parts $\Gamma_\mathrm{in}$, $\Gamma_\mathrm{w}$ and $\Gamma_\mathrm{out}$ with topological conditions which will be concretized in later parts. Boundary conditions on the velocity and normal stress are imposed, that is, for some given $u_v:[0,T]\times\Gamma_\mathrm{in} \to \R^d$,
$u_\sigma:[0,T]\times\Gamma_\mathrm{out} \to \R^d$,
\begin{align}
u_v(t,\xi)&= v(t,\xi) \qquad &&t\in[0,T],\,\xi\in \Gamma_\mathrm{in} \label{eq:oseenc} \\
0&= v(t,\xi) \qquad &&t\in[0,T],\,\xi\in \Gamma_\mathrm{w} \label{eq:oseenc0} \\
u_\sigma(t,\xi)&= \sigma(v,P)(t,\xi)\,\mathrm{n}(\xi) \qquad &&t\in  [0,T],\,\xi\in \Gamma_\mathrm{out},\label{eq:oseend} 
\end{align}
where we set $$\sigma(v,P)(t,\xi) = \sigma(v(t,\xi),P(t,\xi)),\qquad t\in[0,T],\,\xi\in \Omega.
$$
Further, the infinitesimal momentum is initialized at $t=0$, i.e.,
\begin{equation}
p(0,\xi)= p_0(\xi), \qquad\xi\in \Omega,\label{eq:oseene}
\end{equation}
\end{subequations}
where $p_0:\R^3 \to \R^3$ is a given initial infinitesimal momentum. The $\R^{d \times d}$-valued stress tensor is defined by
\begin{align*}
\sigma(v,P) := \mu \left(\nabla v + \nabla v^\top\right) -P\Id_d,
\end{align*}
where $\Id_d\in \R^{d \times d}$ is the identity matrix and $\mu>0$ is the dynamic viscosity of the fluid. 
Since the velocity is divergence-free,
Schwarz's theorem leads to
\begin{align} \label{eq:schwarz}
\div \nabla v^\top = \div \begin{pmatrix}
\nabla v_1\,\ldots\, \nabla v_d
\end{pmatrix}
= \begin{pmatrix}
 \tfrac{\partial}{\partial \xi_1}\div v\\
  \quad\vdots\\
   \tfrac{\partial}{\partial \xi_d}\div v
\end{pmatrix}=0,
\end{align}
whence the stress tensor slightly simplifies to
\begin{align*}
\div \sigma(v,P) = \div( \mu \nabla v - \Id_dP).
\end{align*}
By neglecting the initial and boundary conditions, Oseen's equations read in compact form 
\begin{align}\label{eq:towardsph}
    \begin{pmatrix}
    \dot{p}\\
    0
\end{pmatrix} &= \begin{pmatrix}
        \mu \Delta - \rho b \cdot \nabla & -\nabla\\
        -\div & 0
    \end{pmatrix}
    \begin{pmatrix}
        \frac{1}{\rho} & 0\\
        0&\Id
    \end{pmatrix}
    \begin{pmatrix}
         p\\P
    \end{pmatrix}.
\end{align}
By a~formal consideration, the operator in \eqref{eq:towardsph} is additively composed by skew-adjoint and self-adjoint nonpositive operators, i.e.,
\begin{align*}
    \begin{pmatrix}
        \mu \Delta - \rho b \cdot \nabla & -\nabla\\
        -\div & 0
\end{pmatrix}
=
\begin{pmatrix}
    -\rho b \cdot \nabla & -\nabla\\
    -\div&0
\end{pmatrix}
    + 
    \begin{pmatrix}
        \mu \Delta & 0\\
        0&0
    \end{pmatrix}
\end{align*}
and, moreover, $\frac1\rho p$ is the derivative of the quadratic functional
\begin{align}\label{eq:hamiltonian}
    \mathcal{H}(p) := \int_{\Omega} \frac{1}{2\rho}\|p(\xi)\|^2\,\mathrm{d}\xi
\end{align}
expressing the kinetic energy. Hence, the system \eqref{eq:towardsph} amazingly resembles the structure of port-Hamil\-tonian differential-algebraic equations as considered in \cite{beattie2018linear}. 
In the formulation \eqref{eq:towardsph}, the incompressibility condition can be regarded as an algebraic constraint. 
In our later - strictly analytical and less formal considerations - we will regard Oseen's equation rather as an infinite-dimensional ordinary differential equation. Here, this can be done by including the algebraic constraint in the function space on which the infinitesimal momentum evolves. In the case of models for incompressible flows, this is indeed possible as no control enters directly to the incompressibility condition $\div v=0$.\\
We note that, whereas finite-dimensional port-Hamiltonian DAEs are very well understood, there are many open questions regarding the infinite-dimensional case. Though first promising results have been published recently~\cite{jacob2022solvability,MehrZwar23,reis2021some}, the authors think that this area is still in an embryonic stage, and there is a~need for a~comprehensive theory for infinite-dimensional port-Hamiltonian DAEs. This is however beyond the scope of the present article.

Moreover, we remark that in the typical formulation of Navier-Stokes (and thus also Oseen's) equations, the velocity rather than the infinitesimal momentum is considered as an unknown. We prefer to use the latter since, for mechanical systems in port-Hamiltonian formulation, the state corresponding to kinetic energy is typically given by momenta \cite{JvdS14}. Note that the velocity and momentum in incompressible flow problems only differ by multiplication with the constant mass density of the fluid.

\section{Formal derivation of energy balance: Flow trough a cylinder}\label{eq:cyl}
\begin{figure}[h]
	\begin{tikzpicture}
	\draw (0,0) -- (4,0);
	\draw (0,3) -- (4,3);
	\draw (0,0) arc (270:90:0.5 and 1.5);
	\draw[dashed] (0,0) arc (270:90:-0.5 and 1.5);
	\draw (4,1.5) ellipse (0.5 and 1.5);
	
	\draw[->] (1.0,2.5) -- (2,2.5);
	\draw[->] (2.0,1.8) -- (3,1.8);
	\draw[->] (1.0,1.1) -- (2,1.1);
	
	\draw (2,-.3) node {$\Gamma_\mathrm{w}$};
	\draw (2,3.3) node {$\Gamma_\mathrm{w}$};
	\draw (-.1,1.5) node {$\Gamma_\mathrm{in}$};
	\draw (4.1,1.5) node {$\Gamma_\mathrm{out}$};
	
	\end{tikzpicture}
	\label{fig:cylinder}
	\caption{Cylindrical domain in $\R^3$ with in- and outflow}
\end{figure}
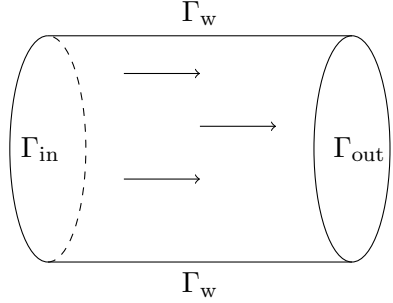

Our motivating example is depicted in Figure~\ref{fig:cylinder} and is given by a tube, i.e., a cylindric domain in three dimensions. Inside the tube, there is a flow field with velocity $v:\R^3\to \R^3$ and pressure $P:\R^3\to \R$. 
 The boundary is partitioned into three distinct parts, namely the side wall and the two tube ends, i.e.,
\begin{align*}
\Gamma = \Gamma_\mathrm{in}\cupdot \Gamma_\mathrm{w} \cupdot \Gamma_\mathrm{out}.
\end{align*}
The left tube end $\Gamma_\mathrm{in}$ is the {inflow boundary},  $\Gamma_\mathrm{w}$ is the side wall, and the right tube end $\Gamma_\mathrm{out}$ is the {outflow boundary}.

\noindent We consider the boundary conditions
\begin{align*}
v &= u_v  &&\text{on } \Gamma_\mathrm{in}  &&\text{prescribed inflow velocity,}\\
v &= 0  &&\text{on } \Gamma_\mathrm{w}  &&\text{no-slip,}\\
\sigma(v,P)\,\mathrm{n} &= u_\sigma &&\text{on } \Gamma_\mathrm{out}  &&\text{prescribed normal stress}.
\end{align*}

In case where $u_\sigma$ vanishes, the latter boundary condition is typically referred to as ``do-nothing condition''.

Next we deduce an energy balance. Note that all our formal derivations can also been done in arbitrary spatial dimensions $d\geq2$, and our analytical results in Section~\ref{sec:oseen_node} indeed cover the general case. In this part, we stick to $d=3$ for illustration purposes.

We first derive an auxiliary result with respect the advection term. As $b$ is divergence free with vanishing normal boundary trace, for a divergence-free function $v:\R^3\to \R^3$ we have, by using the standard inner product in $L^2$,
\begin{align*}
\langle (\rho b\cdot \nabla)v,v\rangle_{L^2(\Omega;\R^3)} &= \rho\sum_{i=1}^3 \langle b\cdot \nabla v_i,v_i\rangle_{L^2(\Omega)} = \rho\sum_{i=1}^3 \langle \nabla v_i,b\,v_i\rangle_{L^2(\Omega;\R^3)}  \\
&= \rho\sum_{i=1}^3 -\langle v_i,  \div (b\, v_i)\rangle_{L^2(\Omega)} + \langle  v_i, \mathrm{n}^\top(b\, v_i)\rangle_{L^2(\Gamma)}\\
&=  \rho\sum_{i=1}^3 -\langle v_i,  b\cdot \nabla v_i)\rangle_{L^2(\Omega)}= - \langle (\rho b\cdot \nabla)v,v\rangle_{L^2(\Omega;\R^3)},
\end{align*}
where the latter equality follows by the chain rule of the divergence $\div(bv_i) = \div(b)v_i + b\cdot \nabla v_i$. The above findings in particular imply that
\begin{align} \label{eq:advection_zero}
 \forall\,v:\,\div v = 0: \quad \langle (\rho b\cdot \nabla)v,v\rangle_{L^2(\Omega;\R^3)} = 0.
\end{align}
Then, by setting $p(t)=p(t,\cdot):\Omega\to\R^3$, and using the  formula \eqref{eq:hamiltonian}
for the kinetic energy $\mathcal{H}(p(t))$, the infinitesimal momentum in \eqref{eq:oseen}, formally fulfills the power balance
\begin{align*}
&\phantom{=}\tfrac{\mathrm{d}}{\mathrm{d}t} \mathcal{H}(p(t)) = \langle \dot p(t), v(t)\rangle\\& = \langle \div \sigma(v(t),P(t)),v(t)\rangle_{L^2(\Omega:\R^3)} + \langle (\rho b\cdot \nabla)v(t),v(t)\rangle_{L^2(\Omega;\R^3)}\\
&=- \langle\sigma(v(t),P(t)), \nabla v\rangle_{L^2(\Omega;\R^{3\times 3})} + \langle \sigma(v(t),P(t))\mathrm{n},v(t)\rangle_{L^2(\Gamma,\R^3)}\\
&= \langle  \mu\nabla v(t) -\Id_3P(t),\nabla v\rangle_{L^2(\Omega,\R^{3\times 3})} + \langle\left(\mu\nabla v(t) -\Id_3P(t)\right)\mathrm{\mathrm{n}},v\rangle_{L^2(\Gamma,\R^3)}.
\end{align*}
Hereby, we note that the inner product $\langle  \cdot,\cdot\rangle_{L^2(\Omega,\R^{d \times d})}$ consists of the integral of the pointwise Hilbert-Schmidt inner product of the two arguments.

A~consequence of the above computations is that the infinitesimal change of energy is additively expressed by a term on the domain and the boundary. Unsurprisingly, each summand has the physical dimension of power: the first one is the volume integral over the product of force per surface element and infinitesimal change of velocity; the second summand is a~surface integral of the  the product of force per surface element and velocity.

Now by using that for any divergence-free $v:\Omega\to \R^3$ and $P:\Omega \to \R$,
\begin{align}\label{eq:Ivnull}
\langle \Id_3P,\nabla v\rangle_{L^2(\Omega;\R^{3\times 3})} = \int_\Omega \tr \left(\begin{pmatrix}P&&\\&P&&\\&&P
\end{pmatrix} \left(\begin{smallmatrix}
\tfrac{\partial v_{1}}{\partial \xi_1}& \tfrac{\partial v_{1}}{\partial \xi_2} & \tfrac{\partial v_{1}}{\partial \xi_3}\\
\tfrac{\partial v_2}{\partial \xi_1}& \tfrac{\partial v_2}{\partial \xi_2} & \tfrac{\partial v_2}{\partial \xi_3}\\
\tfrac{\partial v_3}{\partial \xi_1}& \tfrac{\partial v_3}{\partial \xi_2} & \tfrac{\partial v_3}{\partial \xi_3}\\\end{smallmatrix}\right)\right)\,\mathrm{d}\xi = \int_\Omega P \div v\,\mathrm{d}\xi = 0,
\end{align}	
the above power balance simplifies to
	\begin{align}\label{eq:powerbalance}
\frac{\mathrm{d}}{\mathrm{d}t} \mathcal{H}(p(t)) =  -\mu\|\nabla v(t)\|^2_{L^2(\Omega;\R^{d \times d})} +  \langle  \left(\mu\nabla v(t)-\Id_3P(t)\right)\mathrm{n},v(t)\rangle_{L^2(\Gamma;\R^d)}.
	\end{align}
That is, in the power balance, the pressure only enters on the boundary.
	The term $\mu\|\nabla v(t)\|^2_{L^2(\Omega;\R^{3\times 3})}$ is called \textit{enstrophy} \cite[Chap.~I, Sec.~2]{MariTema1998} and, as can be seen in the above formula \eqref{eq:energybalance}, it is crucial in determining the rate of (interior) kinetic energy decay. 
	By incorporating the split of the boundary into three parts according to Figure~\ref{fig:cylinder}, the boundary term reads (for sake of brevity we omit the time variable)
	\begin{align*}
	&\langle \left(\mu \nabla v-\Id_3 P\right)\mathrm{n},v\rangle_{L^2(\Gamma,\R^3)} \\&=  \langle \left(\mu \nabla v-\Id_3 P\right)\mathrm{n},v\rangle_{L^2(\Gamma_\mathrm{in},\R^3)} + \langle \left(\nu\nabla v -\Id_3 P\right)\mathrm{n},v\rangle_{L^2(\Gamma_\mathrm{w},\R^3)}  + \langle \left(\nu\nabla v -\Id_3 P\right)\mathrm{n},v\rangle_{L^2(\Gamma_\mathrm{out},\R^d)}\\
	&= \langle  \left(\mu\nabla v -\Id_3 P\right)\mathrm{n},u_v\rangle_{L^2(\Gamma_\mathrm{in},\R^3)}
 + \langle  u_\sigma,v\rangle_{L^2(\Gamma_\mathrm{out},\R^3)},
	\end{align*}
	as $v=u_v$ on $\Gamma_\mathrm{in}$ (velocity prescription), $v = 0$ on $\Gamma_\mathrm{w}$ (no-slip condition) and $(\nu \nabla v-\Id_3P)\mathrm{n} = u_\sigma$ on $\Gamma_\mathrm{out}$ (normal stress condition).

Using the latter and integrating \eqref{eq:powerbalance} over $[t_0,t_1]\subset[0,T]$, $t_0\leq t_1$, we obtain the energy balance
	\begin{align}
&\phantom{=}\mathcal{H}(p(t_1))-\mathcal{H}(p(t_0))\nonumber \\&=  -\mu\int_{t_0}^{t_1}\|\nabla v(t)\|^2_{L^2(\Omega;\R^{d \times d})}\mathrm{d}t \nonumber\\&\qquad+ \int_{t_0}^{t_1} \langle  \left(\mu\nabla v(t) -\Id_3 P(t)\right)\mathrm{n},u_v(t)\rangle_{L^2(\Gamma_\mathrm{in},\R^3)}\mathrm{d}t 
 + \int_{t_0}^{t_1}\langle  u_\sigma(t),v(t)\rangle_{L^2(\Gamma_\mathrm{out},\R^3)}\mathrm{d}t \nonumber\\[1mm]&\leq \int_{t_0}^{t_1}\langle  \left(\mu\nabla v(t) -\Id_3 P(t)\right)\mathrm{n},u_v(t)\rangle_{L^2(\Gamma_\mathrm{in},\R^3)}\mathrm{d}t 
 + \int_{t_0}^{t_1}\langle  u_\sigma(t),v(t)\rangle_{L^2(\Gamma_\mathrm{out},\R^3)}\mathrm{d}t .\label{eq:energybalance}
	\end{align}
Energy loss is expressed by the temporal integral over the enstrophy, whereas the energy exchange with outside is formed by the surface integral over the inner product between the velocity and normal component of the stress tensor  
on the boundary.

\begin{rem}
One can straightforwardly include 
 velocity-dependent reaction terms of the form $-c v$ in all our considerations, with a constant $c>0$ in the right-hand side of the dynamics~\eqref{eq:oseena}. This would lead to an additional energy loss term $-c\int_{t_0}^{t_1}\|v\|^2_{L^2(\Omega;\R^d)}\mathrm{d}t$ in the energy balance~\eqref{eq:energybalance}.
\end{rem}

\section{Port-Hamiltonian system nodes}\label{sec:pHnodes}

In this part we present the general analytical background for our approach to port-Hami\-ltonian representation of Oseen's equations.

We first declare some general functional analytic notation. Let $\calX$, $\calY$ be Hilbert spaces. Throughout this article, all Hilbert spaces are real. 
The norm in {$\calX$} will be denoted by $\|\cdot\|_{{\calX}}$ or {simply} $\|\cdot\|$, if clear from context. The identity mapping in $\calX$ is abbreviated by $\Id_{\calX}$ (or
just $\Id$, if clear from context), and we set $\Id_n:=\Id_{\R^n}$.

The symbol $\calX^*$ stands for the {\em dual} of $\calX$; the duality product is denoted by $\<\cdot,\cdot\>_{\calX^*,\calX}$, and the inner product in $\calX$ by $\<\cdot,\cdot\>_{\calX}$. 
We write $L(\calX,\calY)$ for the space of bounded linear operators from $\calX$ to $\calY$. As usual, we abbreviate $L(\calX):= L(\calX,\calX)$. The adjoint of a~densely defined linear operator $A:\calX\supset\dom A\to\calY$ is denoted by $A^*:\calY^*\supset\dom A^*\to\calX^*$.

We adopt the notation of the book \cite{adams2003sobolev} by {\sc Adams} for Lebesgue and Sobolev spaces as well as for the spaces of  continuous and  continuously differentiable functions. For function spaces with values in a~Hilbert space $\calX$, we indicate this by denoting ``$;\calX$'' after specifying the (spatial or temporal) domain. For instance, the Lebesgue space of $p$-integrable $\calX$-valued functions on the domain $\Omega$ is $L^p(\Omega;\calX)$.
Note that, throughout this article, integration of $\calX$-valued functions always has to be understood in the Bochner sense \cite{Diestel77}.

We briefly review the basics of infinite-dimensional port-Hamiltonian systems as presented in \cite{PhilReis23}. In principle, this approach is based on the formulation
\begin{equation}\label{e:phs}
\begin{aligned}
\dot x(t)&=(J\;\:-\;R\:)Hx(t)+(B-P)u(t),\\
y(t)&=(B^\top+P^\top)Hx(t)+(S-N)u(t)
\end{aligned}
\end{equation}
of finite-dimensional port-Hamiltonian systems from \cite{beattie2018linear,MehrMora19}, where
$J\in\R^{n\times n}$ and $N\in\R^{m\times m}$ are skew-symmetric, and $H\in\R^{n\times n}$, $W := \sbmat RP{P^\top}S\in\R^{(n+m)\times (n+m)}$ are symmetric
positive semi-definite.
The total energy of the system is given by the {\em Hamiltonian} $\calH(x) = \frac 12\, x^\top H {x}$. By using \eqref{e:phs} together with skew-symmetry
of $J$ and $N$, as well as positive semi-definiteness of $W$ and $H$, for all $t>0$, $u\in L^2([0,t];\R^m)$, $x_0\in\R^n$, the solution of
\eqref{e:phs} with $x(0)=x_0$ fulfills the {\em dissipation inequality}
\begin{align}
\calH(x(t))-\calH(x_0)=-\int_0^t \spvek {Hx(\tau)}{u(\tau)}^\top\sbmat RP{P^*}S\spvek {Hx(\tau)}{u(\tau)}{\rm d}\tau&+\int_0^t  u(\tau)^\top{y(\tau)}{\rm d}\tau \nonumber \\
&\leq \int_0^t u(\tau)^\top {y(\tau)}{\rm d}\tau,\label{eq:enbalfin}
\end{align}
which has the physical interpretation of an energy balance. Namely, whilst $\calH(x(t))$ stands for the energy stored at time $t$, the first integral after the equality sign is the energy dissipated by the system during the time interval $[0,t]$, whereas $u(\tau)^\top{y(\tau)}$ can be regarded as the external power supply to the system.

\noindent In compact form, \eqref{e:phs} can be rewritten as
\begin{align}\label{e:phs_compact}
\spvek{\dot x(t)}{y(t)} &= \sbmat{\Id_n}{}{}{-\Id_m} M\spvek{Hx(t)}{u(t)},
\end{align}
where $M:=\sbmat{J-R}{B-P}{-B^\top-P^\top}{N-S}\in\R^{(n+m)\times (n+m)}$ is a~dissipative matrix. The generalization to infinite-dimensional systems as presented in \cite{PhilReis23} uses exactly this representation by incorporating the theory of {\em system nodes} by {\sc Staffans} \cite{Staf05}. To allow for partial differential equations, where the dynamics themselves are described by a~differential operator, we clearly presume that $M$ is a~dissipative operator. However,  the incorporation of boundary control requires to allow $M$ to be unbounded and densely defined on some dense subspace of the Cartesian product of the state and input space, whereas the domain is not necessarily a~Cartesian product of a~dense subspace of the state space and the input space. Note that, in \cite{PhilReis23}, also unbounded and non-coercive Hamiltonians are considered. Since the latter is actually not needed for Oseen equations, we ``boil down'' our brief introduction of infinite-dimensional port-Hamiltonian systems to those with bounded and coercive Hamiltonians, which is indeed a~drastic simplification.

The main ingredient for our concept of infinite-dimensional port-Hamiltonian system is given by {\em dissipation nodes}.
To this end, we denote the canonical projection onto $\calX$ and $\calU^*$ in $\calX\times\calU^*$ by $P_{\calX}\in L(\calX\times\calU^*,\calX)$ and $P_{\calU^*}\in L(\calX\times\calU^*,\calU^*)$,
respectively. Further, 
we can canonically identify
$(\calX\times\calU)^*=\calX^*\times\calU^*$. 
\begin{defn}[Dissipation node]\label{def:dissnode}
	Let $\calX$, $\calU$ be Hilbert spaces, let $\calU^*$ be the anti-dual of $\calU$, and identify $\calX$ with its anti-dual, i.e., $\calX^*=\calX$.
	A {\em dissipation node} on $(\calX,\calU)$ is a (possibly unbounded) linear operator $ M:\calX\times\calU\supset\dom M\to\calX\times\calU^*$ that satisfies the following:
	\begin{enumerate}[(a)]
		\item\label{def:dissnode1} $M$ is closed and dissipative;
		\item\label{def:dissnode2} $P_{\calX} M:\calX\times\calU\supset\dom M\to\calX$ is closed;
		\item\label{def:dissnode3} for all $u\in \calU$, there exists some $x'\in \calX$ with $\spvek{x'}{u}\in \dom M$;
		\item\label{def:dissnode4} for the {\em main operator} $F : \calX\supset\dom F\to\calX$ defined by \[\dom F := \setdef{x'\in\calX}{\spvek {x'}0\in\dom  M}\] and $Fx' := P_{\calX}  M\spvek {x'}0$, there exists some $\lambda>0$ such that $\lambda\Id-F$ has dense range.
	\end{enumerate}
\end{defn}
We note that, though the theory and results on port-Hamiltonian systems in \cite{PhilReis23} are formulated and proven for complex Hilbert spaces, the following findings hold as well for the real case, which can be verified by a~complexification.

The previous definition justifies the use of the block operator notation 
\begin{equation}M:=\sbvek{F\&G\\[-1mm]}{K\&L}\label{eq:Mdissnode}\end{equation} for dissipation nodes, where $K\&L=P_{\calX}M$, $F\&G=P_{\calU^*}M$. The symbol ``\&'' is used to indicate that the domain of $M$ is not necessarily a Cartesian product of subspaces of $\calX$ and $\calU$. However, it is shown in \cite{PhilReis23} that, by using the extrapolation space $\calX_{-1}$ defined by the completion of $\calX$ with respect to the norm $\|x\|_{-1}:=\|(\lambda\Id-F)^{-1}\|_{\calX}$ for some $\lambda$ in the resolvent of $F$ (the induced topology is indeed independent on $\lambda$ in the resolvent of $F$), $F$ extends to an operator $F_{-1}\in L(\calX,\calX_{-1})$, and we can regard $G\in L(\calU,\calX_{-1})$. That is,
\begin{equation}F\&G\spvek{x}u=F_{-1}x+Gu\qquad\forall\;\spvek{x}u\in\dom M=\dom (F\&G).\label{eq:FGsplit}\end{equation}
Note that such a~splitting is in general not possible for the operator $K\&L$ defining the output.
\begin{defn}[Port-Hamiltonian system]\label{def:pHnode}
	Let $\calX$, $\calU$ be Hilbert spaces. Assume that $H\in L(\calX)$ is positive, self-adjoint, and has a~bounded inverse. Further, let $M:\calX\times\calU\supset\dom  M\to\calX\times\calU^*$ be a~dissipation node on $(\calX,\calU)$. Then we call
\begin{equation}
\spvek{\dot{x}(t)}{y(t)}
= \sbvek{F\&G\\[-1mm]}{-K\&L} \spvek{{H}{x}(t)}{u(t)}.\label{eq:pHODEnode}\end{equation}
	a~{\em port-Hamiltonian system on $(\calX,\calU)$}.
A {\em classical trajectory} on $[0,T]$ is a triple
\[
(y,x,u)\,\in\, C([0,T];\calU^*)\times C^1([0,T];\calX)\times C([0,T];\calU),
\]
which fulfills \eqref{eq:pHODEnode} pointwise on $[0,T]$, and {\em generalized trajectories} are limits of classical trajectories in the topology of $L^2([0,T];\calY)\times C([0,T];\calX)\times L^2([0,T];\calU)$.
\end{defn}
It follows from \cite[Prop.~3.6]{PhilReis23} that $\sbvek{F\&G\\[-1mm]}{-K\&L}\sbmat{H}00{\Id_{\calU}}$ defines a~system node in the sense of {\sc Staffans} \cite{Staf05}, whence this is called a {\em port-Hamiltonian system node} in \cite{PhilReis23}. Since $H$ is assumed to be bounded, self-adjoint, positive and boundedly invertible, the functional $x\mapsto \|x\|_{H}:=(\langle x(t),Hx(t)\rangle)^{1/2}$ is equivalent to the original norm $\|\cdot\|_\calX$ on $\calX$. It can be moreover seen that $FH$ is a~maximally dissipative (see \cite{EngeNage2000} for a~definition) on $\calX$ endowed with the norm $\|\cdot\|_{H}$ (and thus with the scalar product~$\langle \cdot, H\cdot\rangle$). We collect some properties of port-Hamiltonian systems and their trajectories. Since the result is a~specialization of \cite[Prop.~3.11]{PhilReis23}, we omit the proof.
\begin{prop}[Port-Hamiltonian systems]\label{thm:pHnode}
	Let $\calX$, $\calU$ be Hilbert spaces, and let a~port-Hamil\-tonian system \eqref{eq:pHODEnode} in the sense of Definition~\ref{def:pHnode} on $(\calX,\calU)$ be given. 
 Then the following holds:
	\begin{enumerate}[(a)]
		\item\label{thm:pHnode2} the main operator $A=F{H}$ of ${M}$ generates a~contractive semigroup $\frakA:\R_{\ge0}\to L(\calX)$
in~$(\calX,\|\cdot\|_H)$.
\item\label{thm:pHnode3} For $T{>0}$, $x_0\in \calX$ and $u\in W^{2,1}([0,T];\calU)$ with $\spvek{Hx_0}{u(0)}\in \dom M$, there exists a unique
classical trajectory $(y,x,u)$ for \eqref{eq:pHODEnode} with $x(0)=x_0$.
		\item\label{thm:pHnode4} for all $T>0$, the generalized trajectories $(y,x,u)\,\in\,L^2([0,T];\calU^*)\times C([0,T];\calX)\times  L^2([0,T];\calU)$ (and thus also the classical solutions) fulfill the {\em dissipation inequality}
		\begin{multline}
		\forall\,t\in[0,T]:\qquad \calH(x(t)) - \calH(x(0))\\=
		\int_0^t \left\langle  M\spvek{{H}x(\tau)}{u(\tau)},\spvek {{H}x(\tau)}{u(\tau)}\right\rangle_{\calX\times\calU^*,\calX\times\calU}{\rm d}\tau+\int_0^t \langle {y(\tau)},{u(\tau)}\rangle_{\calU^*,\calU}{\rm d}\tau\\
		\leq\int_0^t \langle {y(\tau)},{u(\tau)}\rangle_{\calU^*,\calU}{\rm d}\tau,\label{eq:enbalinf}
		\end{multline}
		where $\mathcal{H}:\calX\to\R$ with $\calH x := \frac12 \langle x,Hx\rangle$ is the Hamiltonian associated to ${H}$.
	\end{enumerate}
\end{prop}

Though is seems to be exotic in some sense, we emphasize that the input in Proposition~\ref{thm:pHnode}\,\eqref{thm:pHnode3} is truly supposed to be twice weakly differentiable with absolutely integrable second derivative. 
By using that the norms $\|\cdot\|_H$ and $\|\cdot\|_\calX$ are equivalent when $H$ is self-adjoint, bounded and boundedly invertible, 
the semigroup $\frakA:\R_{\ge0}\to L(\calX)$ is bounded when $\calX$ is equipped with the standard norm $\|\cdot\|_\calX$.
We further give some comments on the solution concept, which follow from the results in \cite{PhilReis23}.

\begin{rem}[Classical/generalized trajectories]\label{rem:sols}
Let $T{>0}$, let $\calX$, $\calU$ be Hilbert spaces, and let a~port-Hamilto\-nian system on $(\calX,\calU)$ in the sense of Definition~\ref{def:pHnode} be given.
\begin{enumerate}[(a)]
\item If $(x,u,y)$ is a~classical trajectory, then
\[
\spvek xu\in C([0,T];\dom M),
\]
where $\dom M$ is equipped with the graph norm of $M$.
\item\label{rem:sols1} 
For the split of $F\&G$ as in~\eqref{eq:FGsplit}, consider $F_{-1}H:\calX_{-1}\supset \calX\to \calX_{-1}$. Then $F_{-1}H$ generates a~semigroup $\frakA:\R_{\ge0}\to L(\calX_{-1})$ on $\calX_{-1}$. 
Then $(x,u,y)$ is a~generalized trajectory for \eqref{eq:pHODEnode}, if $x\in C([0,T];\calX)$ and
\begin{equation}\label{eq:mildsol}
\forall\,t\in[0,T]:\quad x(t)=\frakA(t)x(0)+\int_0^t \frakA_{-1}(t-\tau)Gu(\tau){\rm d}\tau,
\end{equation}
where the latter summand has to be interpreted as an integral in the space $\calX_{-1}$. 
\item\label{rem:sols3} If $(x,u,y)$ is a~generalized trajectory for \eqref{eq:pHODEnode}, then \eqref{eq:mildsol} holds. The output evaluation $y(t)=K\&L \spvek {Hx(t)}{u(t)}$ is -- at a glance -- not necessarily well-defined for all $t\in[0,T]$. However, it is shown in \cite[Lem.~4.7.9]{Staf05} (for the larger context of system nodes) that the second integral of $\spvek {x}{u}$ is continuous as a~mapping from $[0,T]$ to $\dom(K\&L)=\dom M$. As a~consequence, the output can -- in the distributional sense -- be defined as the second derivative of $K\&L$ applied to the second integral of $\spvek {Hx}{u}$. This can be used to show that $(x,u,y)$ is a~generalized trajectory for \eqref{eq:pHODEnode} if, and only if, \eqref{eq:mildsol} and
\[
y=\left(t\mapsto\tfrac{{\rm d}^2}{{\rm d}t^2}\,K\&L\int_0^t(t-\tau)\spvek {Hx(\tau)}{u(\tau)}{\rm d}\tau\right)\in L^2([0,T];\calY).
\]
\end{enumerate}
\end{rem}

\section{Oseen equations as port-Hamiltonian system}
\label{sec:oseen_node}
In this part, we formulate the PDE \eqref{eq:oseen} by means of a port-Hamiltonian system in the sense of Definition~\ref{def:pHnode}.  In view of the representation \eqref{eq:towardsph} and the quadratic Hamiltonian energy as in~\eqref{eq:hamiltonian}, 
we define the Hamiltonian operator given as a~simple multiplication operator
\begin{align*}
    Hz = \tfrac{1}{\rho} z.
\end{align*}
such that $\mathcal{H}(p) = \frac12\langle p,Hp\rangle_{L^2(\Omega;\R^d)}$.
Since the density $\rho>0$ is assumed to be constant, this operator boundedly maps square integrable and divergence-free momenta to square integrable and divergence-free velocities. 
 Having defined the Hamiltonian, it remains to describe the dissipation node incorporating the differential operators in the formulation \eqref{eq:towardsph} to obtain a port-Hamiltonian system in the form of \eqref{eq:pHODEnode}. This will be done in the sequel.

\subsection{Assumptions on the domain}\label{sec:domain}

The following assumptions are made on the domain $\Omega$ and (parts of the) boundary $\Gamma,\Gamma_\mathrm{in},\Gamma_\mathrm{w},\Gamma_\mathrm{out}$ throughout the remainder of this work: 
\begin{itemize}
    \item $\Omega\subset \mathbb{R}^d$, $d\ge2$ is a~Lipschitz domain with boundary $\Gamma$. 
    \item $\Gamma_\mathrm{in},\Gamma_\mathrm{w},\Gamma_\mathrm{out}\subset\Gamma$ are relatively open and Lipschitz domains of dimension $d-1$ with boundary.
    \item $\Gamma_\mathrm{in}$ and $\Gamma_\mathrm{out}$ are non-empty.
    \item $\Gamma_\mathrm{in}$, $\Gamma_\mathrm{w}$, $\Gamma_\mathrm{out}$ are disjoint, that is, 
    \[\Gamma_\mathrm{in}\cap\Gamma_\mathrm{w}=\Gamma_\mathrm{in}\cap\Gamma_\mathrm{out}=\Gamma_\mathrm{w}\cap\Gamma_\mathrm{out}=\emptyset.\]
    \item The closures of $\Gamma_\mathrm{in}$, $\Gamma_\mathrm{w}$, $\Gamma_\mathrm{out}$ unify to the whole boundary, i.e., \[\overline{\Gamma_\mathrm{in}}\cup\overline{\Gamma_\mathrm{w}} \cup\overline{\Gamma_\mathrm{out}}=\Gamma.\]
    \item $\Gamma_\mathrm{in}$ and $\Gamma_\mathrm{out}$ do not touch, that is,
\[\overline{\Gamma_\mathrm{in}}\cap \overline{\Gamma_\mathrm{out}} = \emptyset.\]
\end{itemize}

In view of the tube depicted in Figure~\ref{fig:cylinder}, $\Gamma_\mathrm{w}$ is the wall of the cylinder as the no-slip boundary with vanishing velocities and $\Gamma_\mathrm{in}$ is the inflow boundary with prescribed velocities. Likewise, $\Gamma_\mathrm{out}$ will be the part of the boundary in which conditions on the stress tensor $\sigma(v,P)$ in normal direction are imposed.  In the example illustrated in Figure~\ref{fig:cylinder}, this is the part containing the ``do-nothing conditions''. In fact, we also allow nonzero conditions on $\Gamma_\mathrm{out}$ in our forthcoming analysis.

\subsection{Spaces and traces}

Let $H^{1/2}({\Gamma};\R^d)$ be the Sobolev space of fractional order 1/2. Further, for some relatively open subset $\Gamma_i\subset\Gamma$, the space $H^{1/2}_0({\Gamma_i};\R^d)$ consists of elements of $H^{1/2}({\Gamma_i};\R^d)$ whose extension by zero is an element of $H^{1/2}({\Gamma};\R^d)$ \cite{adams2003sobolev}.

We consider the {\em trace operator}
\begin{align}\label{eq:dirichlet_fullboundary_whole}
\gamma: H^1(\Omega;\R^d) \to H^{1/2}( \Gamma;\R^d),
\end{align}
which is bounded and surjective by \cite[Thm.~1.5.1.3]{Gris85}. 
Further, consider
the space
\begin{align*}
H^1_{\sigma,\mathrm{w}}(\Omega;\R^d) &=  \{v\in H^1(\Omega;\R^d) \,|\, \div v = 0\ \mathrm{ on }\ \Omega\text{ and } \gamma v\vert_{{\Gamma}_\mathrm{w}}=0\}.
\end{align*}
It can be seen that the operator defining the restriction $w\vert_{{\Gamma}_\mathrm{w}}\in H^{1/2}({\Gamma_\mathrm{w}};\R^d)$ of $w\in H^{1/2}(\Gamma;\R^d)$ is bounded. This together with boundedness of the trace operator $\gamma$  and boundedness of the divergence operator $\div : H^1(\Omega;\R^d) \to L^2(\Omega;\R)$ yields that $H^1_{\sigma,\mathrm{w}}(\Omega;\R^d)$ is a closed subspace of $H^1(\Omega;\R^d)$; it is therefore a~Hilbert space endowed with the inner product in $H^1(\Omega;\R^d)$.

Since, by our assumptions, $\Gamma_\mathrm{in}$ and $\Gamma_\mathrm{out}$ do not touch, there is a canonical isomorphism
\begin{equation}\left\{w\in H^{1/2}(\Gamma;\R^d)\,|\,w\vert_{{\Gamma}_\mathrm{w}}=0\right\}\cong H_0^{1/2}(\Gamma_\mathrm{in};\R^d)\times H_0^{1/2}(\Gamma_\mathrm{out};\R^d).\label{eq:caniso}\end{equation}
Now Gau\ss' divergence theorem yields that $\gamma$ restricts to
\begin{multline}
\label{eq:dirichlet_fullboundary}
{\gamma}: H^1_{\sigma,\mathrm{w}}(\Omega;\R^d) \to\\ 
\left\{(w_\mathrm{in},w_\mathrm{out})\in H_0^{1/2}(\Gamma_\mathrm{in};\R^d)\times H_0^{1/2}(\Gamma_\mathrm{out};\R^d)\,|\,\int_{{\Gamma_\mathrm{in}}}w_\mathrm{in}\cdot \mathrm{n}\,\mathrm{d}\xi+\int_{{\Gamma_\mathrm{out}}} w_\mathrm{out}\cdot \mathrm{n}\,\mathrm{d}\xi=0 \right\}, 
\end{multline}
which is again denoted by $\gamma$ for sake of brevity. Note that the target space of the latter operator is a~closed subspace of $H_0^{1/2}(\Gamma_\mathrm{in};\R^d)\times H_0^{1/2}(\Gamma_\mathrm{out};\R^d)$ with co-dimension one. It is thus a~Hilbert space equipped with the inner product in $H_0^{1/2}(\Gamma_\mathrm{in};\R^d)\times H_0^{1/2}(\Gamma_\mathrm{out};\R^d)$. Indeed, it can be concluded from \cite[Lem.~2.2]{GiraRavi12} that ${\gamma}$ as in \eqref{eq:dirichlet_fullboundary} is surjective.\\
Now we denote 
\begin{align}\label{eq:rest_dir_trace}
\begin{split}
        \gamma_\mathrm{in}:\quad H^1_{\sigma,\mathrm{w}}(\Omega;\R^d)&\to H_0^{1/2}(\Gamma_\mathrm{in};\R^d),\\
    \gamma_\mathrm{out}:\quad H^1_{\sigma,\mathrm{w}}(\Omega;\R^d)&\to H_0^{1/2}(\Gamma_\mathrm{out};\R^d)
    \end{split}
\end{align}
to be the canonical projections of $\gamma$ onto $H_0^{1/2}(\Gamma_\mathrm{in};\R^d)$ and $H_0^{1/2}(\Gamma_\mathrm{out};\R^d)$, respectively. Since, by our assumption on the domain, $\Gamma_\mathrm{in}$ and $\Gamma_\mathrm{out}$ are both non-empty, it follows from the surjectivity of $\gamma$ in \eqref{eq:dirichlet_fullboundary} that both $\gamma_\mathrm{in}$ and $\gamma_\mathrm{out}$ are surjective.

Moreover, let $H(\div,\Omega;\R^{d\times d})$ be the space of square integrable functions with values in $\R^{d \times d}$ for which the row-wise divergence exists and is square integrable, where the divergence $\div$ is defined in a weak sense. That is, for the space $H^1_0(\Omega;\R^d)$ of all elements of $H^1(\Omega;\R^d)$ with vanishing boundary trace,
\[z=\div x\quad\Longleftrightarrow\quad \forall \varphi\in
H^1_0(\Omega;\R^d):\;-\langle\nabla\varphi,x\rangle_{L^2(\Omega;\R^{d \times d})}=\langle\varphi,z\rangle_{L^2(\Omega;\R^d)}.\]
Likewise, $H(\div,\Omega;\R^{d})$ denotes the space of square integrable $\R^d$-valued functions with square integrable divergence, where the latter is defined via
\[z=\div x\quad\Longleftrightarrow\quad \forall \varphi\in
H^1_0(\Omega):\;-\langle\nabla\varphi,x\rangle_{L^2(\Omega;\R^{d})}=\langle\varphi,z\rangle_{L^2(\Omega)}.\]

Here we collect the fundamentals to define the normal trace of elements of $H(\div,\Omega;\R^{d\times d})$. 
First - according to the notation in \cite{adams2003sobolev} - we denote, for $\hat\Gamma \in \{\Gamma,\Gamma_\mathrm{in}, \Gamma_\mathrm{out}\}$,
 \begin{align*}
     H^{-1/2}(\hat{\Gamma};\R^d)&=H^{1/2}_0(\hat{\Gamma};\R^d)^*.
 \end{align*}
Since the boundary of $\Gamma$ is empty, we have $H^{1/2}_0({\Gamma};\R^d)=H^{1/2}({\Gamma};\R^d)$.\\
Using surjectivity of the trace operator $\gamma$ from $H^{1}(\Omega;\R^d)$ to $H^{1/2}(\Gamma;\R^d)$, the \textit{normal trace} 
$\gamma_{\mathrm{n}} \sigma\in H^{-1/2}(\Gamma;\R^d)$
of $\sigma \in H(\div,\Omega;\R^{d\times d})$ is well-defined by 
\[\gamma_{\mathrm{n}} \sigma\in H^{-1/2}(\Gamma;\R^d),\]
where
\begin{align}\label{eq:normaltrace}
\begin{split}
\forall z\in H^1(\Omega;\R^d)\,:\,\langle \gamma_{\mathrm{n}}\sigma,\gamma z &\rangle_{H^{-1/2}(\Gamma;\R^d),H^{1/2}(\Gamma;\R^d)}\\&= \langle \div \sigma,z\rangle_{L^2(\Omega;\R^d)} + \langle \sigma,\nabla z\rangle_{L^2(\Omega;\R^{d\times d})}.
\end{split}
\end{align}
In the case where $\sigma$ and $\Gamma$ are smooth, Gau\ss' divergence theorem yields that $w(\xi) = \mathrm{n}(\xi)^\top v(\xi)$  for all $\xi\in \Gamma$, cf.\ \cite[Chap.~16]{Tart07}.\\
The normal traces of $\sigma \in H(\div,\Omega;\R^{d\times d})$ on ${\Gamma}_\mathrm{in}$ and ${\Gamma}_\mathrm{out}$, which are denoted by
\[\gamma_{\mathrm{n},\mathrm{in}} \sigma\in H^{-1/2}({\Gamma}_\mathrm{in};\R^d),\quad \gamma_{\mathrm{n},\mathrm{out}} \sigma\in H^{-1/2}({\Gamma}_\mathrm{out};\R^d),\]
 can be be defined via
 testing with weakly differentiable functions vanishing outside 
${\Gamma}_\mathrm{in}$ and ${\Gamma}_\mathrm{out}$, respectively. That is,
\begin{align*}
\forall z\in H^1(\Omega;\R^d)\text{ with }\gamma z\vert_{\Gamma_\mathrm{w}\cup \Gamma_\mathrm{out}}=0\,:\qquad\qquad\\
\langle \gamma_{\mathrm{n},\mathrm{in}}\sigma,\gamma_{\mathrm{in}} z \rangle_{H^{-1/2}({\Gamma}_\mathrm{in};\R^d),H^{1/2}({\Gamma}_\mathrm{in};\R^d)}&= \langle \div \sigma,z\rangle_{L^2(\Omega;\R^d)} + \langle \sigma,\nabla z\rangle_{L^2(\Omega;\R^{d\times d})},\\[1mm]
\forall z\in H^1(\Omega;\R^d)\text{ with }\gamma z\vert_{\Gamma_\mathrm{w}\cup \Gamma_\mathrm{in}}=0\,:\qquad\qquad\\\langle \gamma_{\mathrm{n},\mathrm{out}}\sigma,\gamma_{\mathrm{out}} z \rangle_{H^{-1/2}({\Gamma}_\mathrm{out};\R^d),H^{1/2}({\Gamma}_\mathrm{out};\R^d)}&= \langle \div \sigma,z\rangle_{L^2(\Omega;\R^d)} + \langle \sigma,\nabla z\rangle_{L^2(\Omega;\R^{d\times d})}.
\end{align*}
Note that in the above definition, $\gamma_\mathrm{in}$ and $\gamma_\mathrm{out}$ denote the surrogates of the restricted trace as defined in \eqref{eq:rest_dir_trace} for functions $H^1(\Omega;\R^3)$, mapping surjectively onto $H^{1/2}_0(\Gamma_\mathrm{in};\R^3)$ and $H^{1/2}_0(\Gamma_\mathrm{out};\R^3)$ respectively.\\
A~consequence is the following Green's identity, which will play a central role in the remainder:
\begin{align}\label{eq:normaltracehat}
\begin{split}
&\!\!\!\!\forall z\in H^1_{\sigma,\mathrm{w}}(\Omega;\R^d),\, \sigma\in H(\div,\Omega;\R^{d\times d})\,:\qquad\qquad\\
&\phantom{=}\langle \gamma_{\mathrm{n},\mathrm{in}}\sigma,\gamma_{\mathrm{in}} z \rangle_{H^{-1/2}({\Gamma}_\mathrm{in};\R^d),H^{1/2}({\Gamma}_\mathrm{in};\R^d)}
+\langle \gamma_{\mathrm{n},\mathrm{out}}\sigma,\gamma_{\mathrm{out}} z \rangle_{H^{-1/2}({\Gamma}_\mathrm{out};\R^d),H^{1/2}({\Gamma}_\mathrm{out};\R^d)}
\\&= \langle \div \sigma,z\rangle_{L^2(\Omega;\R^d)} + \langle \sigma,\nabla z\rangle_{L^2(\Omega;\R^{d\times d})}.
\end{split}
\end{align}

\subsection{Dissipation node corresponding to the Oseen system}
We now introduce the system node corresponding to the Oseen problem \eqref{eq:oseen}. 
Hereby, we assume that $\Omega\subset\R^d$, $d\ge2$, with boundary $\Gamma$, and boundary parts $\Gamma_\mathrm{in},\Gamma_\mathrm{w},\Gamma_\mathrm{out}$ having the properties defined in Section~\ref{sec:domain}. Further, we assume that $b\in L^\infty(\Omega;\R^d)\cap H(\div,\Omega;\R^d)$ with $\div b=0$ and trivial normal trace, i.e., $\gamma_{\mathrm{n}}b=0$. Further, $\mu,\rho\in\R$ are positive constants for the dynamic viscosity and the mass density.

The remainder of this section is devoted to prove that the following is a~dissipation node on $(\calX,\calU) = \left(L^2_\sigma(\Omega;\R^3), H^{1/2}_{0}(\Gamma_{\mathrm{in}};\R^d)\times H^{-1/2}_{0}(\Gamma_\mathrm{out};\R^d)\right)$ with the state space $L^2_\sigma(\Omega;\R^d) = \{v\in H(\div,\Omega;\R^d)\,|\,\div v = 0\}$:
\begin{subequations}\label{eq:dissnode_oseen}
	\begin{align}\begin{split}
	\dom M\!=\!\dom{F\&G}\!:=\!\bigg\{&\spvekk{v}{u_v}{u_\sigma}\in H^1_{\sigma,{\mathrm{w}}}(\Omega;\R^d)\times H^{1/2}_{0}(\Gamma_{\mathrm{in}};\R^d)\times H^{-1/2}_{0}(\Gamma_\mathrm{out};\R^d) \, \bigg\vert  \\ &  \qquad \exists P\in L^2(\Omega)\ \text{s.t.}\ \mu \nabla v - P\Id_d \in H(\div,\Omega;\R^{d\times d})\;\\&  \qquad \wedge\;\gamma_{\mathrm{in}}v=u_v\;\wedge\;\gamma_{\mathrm{n},\mathrm{out}}(\mu \nabla v-P\Id_d)=u_\sigma\bigg\}
	\end{split}
	\end{align}
	and
	\begin{align}
	\forall \, \spvekk{v}{u_v}{u_\sigma}\in\dom M:\quad F\&G\spvekk{v}{u_v}{u_\sigma}&=\div (\mu \nabla v  - P\Id_d) + \rho(b\cdot\nabla)v,\label{eq:FG}\\
	K\&L\spvekk{v}{u_v}{u_\sigma}&=-\spvek{\gamma_{\mathrm{n},\mathrm{in}} (\mu \nabla v-P\Id_d)}{\gamma_{\mathrm{out}} v},\;\;
	M=\sbvek{F\&G}{K\&L}.
	\end{align}
\end{subequations}

Crucial in showing that \eqref{eq:dissnode_oseen} is a dissipation node is the following simplified and weak version of De Rham's theorem characterizing the orthogonal complement of divergence free functions as gradient fields, cf.\ \cite[Thm.~2.3]{GiraRavi12} or \cite[Rem.~1.9]{Tema01}. This result ensures well-definition of $F\&G$ as it implies uniqueness of the pressure gradient in the definition \eqref{eq:FG}, cf.\ also \cite{ArenElst12}.
\begin{lem}\label{lem:pressure_recovery}
Let $\Omega\subset\R^d$ be a~bounded Lipschitz domain, and let $f\in H^{-1}(\Omega;\R^d) = H^1_0(\Omega;\R^d)^*$, such that $\langle f,v\rangle_{L^2(\Omega;\R^d)}=0$ for all $v\in H^1_0(\Omega;\R^d)$ with $\div v=0$. Then there exists some $P\in L^2(\Omega)$ such that $f = \nabla P$ holds in the distributional sense. That is, 
\[\forall\, z\in H^1_0(\Omega;\R^d):\quad\langle f,z\rangle_{H^{-1}(\Omega;\R^d),H^1_0(\Omega;\R^d)}
=-\langle P,\div z\rangle_{L^{2}(\Omega)}.
\]
Further, $P$ is unique up to addition of a~constant function.
\end{lem}
We further advance a~lemma which basically states that the drift term $\rho(b\cdot\nabla)v$ does not contribute to dissipation. It is the weak form of eq.\ \eqref{eq:advection_zero}.
\begin{lem}\label{lem:drift}
Let $\Omega\subset\R^d$ be a~bounded Lipschitz domain, and 
let $b\in L^\infty(\Omega)\cap H(\div,\Omega)$ with $\div b=0$ and trivial normal trace, that is, $\gamma_{\mathrm{n}}b=0$. Then for all $v\in H^1(\Omega;\R^d)$,
\[\langle \,v,\rho(b\cdot \nabla) v\rangle_{L^2(\Omega;\R^d)}=0.\]
\end{lem}
\begin{proof}
By boundedness of $\Omega$ and the fact that smooth functions are contained in both $L^1(\Omega;\R^d)$ and $L^2(\Omega;\R^d)$, we have that $L^2(\Omega;\R^d)$ is dense in $L^1(\Omega;\R^d)$.
  As a~consequence, our condition on $b$ is equivalent to
\begin{equation}b\in L^{\infty}(\Omega;\R^d)\;\;\wedge\;\; \forall\,z\in W^{1,1}(\Omega;\R^d):\; {\langle\nabla z,b\rangle_{L^1(\Omega;\R^d),L^\infty(\Omega;\R^d)}}=0,\label{eq:divorth}\end{equation}
where the latter stands for the canonical duality product of $L^1$ and $L^\infty\cong (L^1)^*$.
Let $v\in H^1(\Omega;\R^d)$. The product rule for weak derivatives \cite[Thm.~4.25]{Alt16} yields that the pointwise inner product $v\cdot v\in W^{1,1}(\Omega)$ with
\[\tfrac12\nabla( v\cdot v)=v\cdot \nabla{v}\in L^1(\Omega;\R^{d}),\]
and thus
\begin{align*}
\langle \,v,(b\cdot \nabla) v\rangle_{L^2(\Omega;\R^d)}
&=\langle \,b,v\cdot \nabla v\rangle_{L^2(\Omega;\R^d)}\\
&=\tfrac12\langle \,b, \nabla( v\cdot v)\rangle_{L^\infty(\Omega;\R^d),L^1(\Omega;\R^d)}{\stackrel{\eqref{eq:divorth}}{=}}0.
\end{align*}
\end{proof}
Now we are able to show that $M$ as in \eqref{eq:dissnode_oseen} is a~dissipation node.

\begin{thm}[Dissipation node for Oseen equations]
	Let $\Omega\subset\R^d$, $d\ge2$, with boundary $\Gamma$, and $\Gamma_\mathrm{in},\Gamma_\mathrm{w},\Gamma_\mathrm{out}\subset \Gamma$ satisfying the properties as constituted in Section~\ref{sec:domain}. Further, let $b\in L^\infty(\Omega;\R^d)\cap H(\div,\Omega;\R^d)$ with $\div b=0$ and $\gamma_{\mathrm{n}}b=0$, and let $\mu,\rho\in\R_{>0}$.\\
 Then the operator $M$ as in \eqref{eq:dissnode_oseen} is a~dissipation node on $(\calX,\calU)$ with $\calX = L^2(\Omega;\R^d)$ and $\calU=H^{1/2}_0(\Gamma_\mathrm{in};\R^d)\times H^{-1/2}(\Gamma_\mathrm{out};\R^d)$.
\end{thm}
\begin{proof}
	We successively show that $M$ fulfills \eqref{def:dissnode1}-\eqref{def:dissnode4} in Definition~\ref{def:dissnode}.
	\begin{itemize}
		\item[\eqref{def:dissnode1}] {\em Step~1:} Let $\spvekk{v}{u_v}{u_\sigma}\in\dom M$. 
  We first note that - for obvious reasons - \eqref{eq:Ivnull} holds for weakly differentiable and divergence-free $\R^d$-valued functions as well, whence 
  \begin{equation}\forall\, w\in H^1_{\sigma,{\mathrm{w}}}(\Omega;\R^d):\quad\langle \nabla w,P\Id_d\rangle_{L^2(\Omega;\R^{d\times d})}=0.\label{eq:wdiv0}\end{equation}
Now dissipativity of $M$ follows from (for sake of brevity we neglect the subindices indicating the spaces specifying the duality product)
\begin{align*}
		&\phantom{=}\left\langle \spvekk{v}{u_v}{u_\sigma},M \spvekk{v}{u_v}{u_\sigma}\right\rangle\\
		&{=}\left\langle \spvekk{v}{u_v}{u_\sigma},\spvekk{\div (\mu \nabla - P\Id_d) + \rho(b\cdot\nabla)v}{-\mu \gamma_{\mathrm{n},\mathrm{in}}(\mu \nabla v-P\Id_d)}{-\gamma_{\Gamma_\mathrm{out}}v}
		\right\rangle\\
		&{=}\langle v, \div (\mu \nabla v - P\Id_d)\rangle+\underbrace{\langle v,\rho(b\cdot \nabla) v\rangle}_{{=}0\text{ by Lemma~\ref{lem:drift}}}
		-\langle u_v, \mu \gamma_{\mathrm{n},\mathrm{in}} (\mu\nabla v-P\Id_d)\rangle - \langle u_\sigma, \gamma_{\mathrm{n},\mathrm{out}}v \rangle\\
  		&{=}\langle v, \div (\mu \nabla v - P\Id_d)\rangle
		-\langle \gamma_{\mathrm{in}}v,  \gamma_{\mathrm{n},\mathrm{in}} (\mu\nabla v-P\Id_d)\rangle - \langle \gamma_{\mathrm{out}}v, \gamma_{\mathrm{n},\mathrm{out}}(\mu \nabla v - P\Id_d) \rangle\\
		&{\stackrel{\eqref{eq:normaltracehat}}{=}}-\mu \langle \nabla v, \mu\nabla v-P\Id_d\rangle\\
		&{\stackrel{\eqref{eq:wdiv0}}{=}}- \mu \langle \nabla v, \nabla v\rangle \leq 0.
		\end{align*}
		{\em Step~2:} We show that $ M-\calR$ is onto, where $\calR:\calX\times\calU\to\calX\times\calU^*$ and in a~matrix notation, $\calR=\mathrm{diag}(\Id_\calX,\calR_\mathrm{in},\calR_\mathrm{out}^{-1})$ with the Riesz isomorphisms
\[  \begin{array}{rrcl}
  \calR_\mathrm{in}:& H^{1/2}_0(\Gamma_\mathrm{in};\R^d)&\to &H^{-1/2}(\Gamma_\mathrm{in};\R^d),\\
    \calR_\mathrm{out}:& H^{1/2}_0(\Gamma_\mathrm{out};\R^d)&\to& H^{-1/2}(\Gamma_\mathrm{out};\R^d).
  \end{array}\]
  Let $z\in \calX=L^2(\Omega;\R^d)$ and $(w_1,w_2)\in H^{-1/2}(\Gamma_\mathrm{in};\R^d)\times H^{1/2}_0(\Gamma_\mathrm{out};\R^d)$. We have to find $(x,u_v,u_\sigma)\in\dom M$ such that
		\begin{equation}\big( M-\calR\big)\spvekk{v}{u_v}{u_\sigma}=\spvekk{z}{w_1}{w_2}.\label{eq:maxdiss}\end{equation}
Defining the bilinear form and the linear functional 
\begin{align*}
		  b(v,\varphi) &:= \mu \langle \nabla\varphi,\nabla v\rangle_{L^2(\Omega;\R^{d \times d})}
		+ \langle \varphi, v\rangle_{L^2(\Omega;\R^d)}-\langle \varphi, \rho(b\cdot \nabla) v\rangle_{L^2(\Omega;\R^d)}\\
		&\qquad \qquad +\langle \gamma_{\mathrm{in}}\varphi, \gamma_{\mathrm{in}}v\rangle_{H^{1/2}_0(\Gamma_\mathrm{in};\R^d)} +\langle \gamma_{\mathrm{out}}\varphi, \gamma_{\mathrm{out}}v \rangle_{H_0^{1/2}(\Gamma_\mathrm{out};\R^d)},\\
		F_{z,w_1,w_2}(\varphi)&  :=-\langle \varphi,z\rangle_{L^2(\Omega;\R^d)}-\langle \gamma_{\mathrm{in}}\varphi,w_1\rangle_{H_0^{1/2}(\Gamma_\mathrm{in};\R^d), H^{-1/2}(\Gamma_\mathrm{in};\R^d)}\\&\qquad \qquad -\langle \gamma_{\mathrm{out}}\varphi,w_2\rangle_{H_0^{1/2}(\Gamma_\mathrm{out};\R^d)},
  \end{align*}
we can conclude from the Lax-Milgram lemma \cite[Lem.~2.2.1.1]{Gris85} that there exists some $v\in H_{\sigma,\mathrm{w}}^1(\Omega;\R^d)$, such that,
		\begin{align}
        b(v,\varphi) = F_{z,w_1,w_2}(\varphi) \label{eq:laxmil1}
		\end{align}
   for all $\varphi\in H_{\sigma,\mathrm{w}}^1(\Omega;\R^d)$.
In particular, for all test functions $\varphi \in H^1_0(\Omega;\R^d)$ with $\div\varphi=0$, we have 
$$\langle \div \mu \nabla v -v+ \rho(b\cdot \nabla)v-z,\varphi\rangle_{H^{-1}(\Omega;\R^d), H^1_{0}(\Omega;\R^d) } = 0.
$$ 
An application of Lemma~\ref{lem:pressure_recovery} gives rise to the existence of some $P\in L^2(\Omega)$ with $\div \mu \nabla v - v + \rho(b\cdot \nabla)v- z = \nabla P$. This means that $\div(\mu \nabla v - P\Id_d) = v-\rho(b\cdot \nabla v) + z \in L^2(\Omega)$, i.e., \begin{equation}\div (\mu \nabla v -P\Id_d) + \rho(b\cdot \nabla)v -v= z.\label{eq:sysnodeeq1}\end{equation}
By abbreviating  $\sigma:=\mu \nabla v -P\Id_d$ and by setting $u_v = \gamma_{\mathrm{in}}v$ and 
$u_\sigma= \gamma_{\mathrm{n},{\mathrm{out}}}\sigma$, we compute
\begin{align*}
    &\phantom{=} \langle \gamma_{\mathrm{in}}\varphi, u_v\rangle_{H^{1/2}_0(\Gamma_\mathrm{in};\R^d)} +\langle \gamma_{\mathrm{out}}\varphi, \gamma_{\mathrm{out}}v \rangle_{H_0^{1/2}(\Gamma_\mathrm{out};\R^d)}\\
    &\ \ \stackrel{\eqref{eq:wdiv0}}{=}  \mu \langle \nabla\varphi,\nabla v\rangle_{L^2(\Omega;\R^{d \times d})}-\langle \nabla\varphi,\sigma\rangle_{L^2(\Omega;\R^{d \times d})}\\
    &\qquad  \qquad+\langle \gamma_{\mathrm{in}}\varphi, u_v\rangle_{H^{1/2}_0(\Gamma_\mathrm{in};\R^d)} +\langle \gamma_{\mathrm{out}}\varphi, \gamma_{\mathrm{out}}v \rangle_{H_0^{1/2}(\Gamma_\mathrm{out};\R^d)}\\
    &\underset{\&\eqref{eq:sysnodeeq1}}{\overset{\eqref{eq:laxmil1}}{=}}  
    -\langle \varphi,\div \sigma + \rho(b\cdot \nabla)v -v\rangle_{L^2(\Omega;\R^d)}-\langle \nabla\varphi,\sigma\rangle_{L^2(\Omega;\R^{d \times d})}\\
    &\qquad \qquad-\langle \varphi, v\rangle_{L^2(\Omega;\R^d)}+\langle \varphi, \rho(b\cdot \nabla) v\rangle_{L^2(\Omega;\R^d)}\\
    &\qquad \qquad-\langle \gamma_{\mathrm{in}}\varphi,w_1\rangle_{H_0^{1/2}(\Gamma_\mathrm{in};\R^d), H^{-1/2}(\Gamma_\mathrm{in};\R^d)}-\langle \gamma_{\mathrm{out}}\varphi,w_2\rangle_{H_0^{1/2}(\Gamma_\mathrm{out};\R^d)}\\
	& \quad  =-\langle \varphi,\div \sigma\rangle_{L^2(\Omega;\R^d)}-\langle \nabla\varphi,\sigma\rangle_{L^2(\Omega;\R^{d \times d})}\\
    &\qquad\qquad-\langle \gamma_{\mathrm{in}}\varphi,w_1\rangle_{H_0^{1/2}(\Gamma_\mathrm{in};\R^d), H^{-1/2}(\Gamma_\mathrm{in};\R^d)}-\langle \gamma_{\mathrm{out}}\varphi,w_2\rangle_{H_0^{1/2}(\Gamma_\mathrm{out};\R^d)}\\
    &\ \ \stackrel{\eqref{eq:normaltracehat}}{=}-
    \langle \gamma_{\mathrm{n},\mathrm{in}}\sigma,\gamma_{\mathrm{in}} \varphi \rangle_{H^{-1/2}({\Gamma}_\mathrm{in};\R^d),H^{1/2}({\Gamma}_\mathrm{in};\R^d)}
    -\langle u_\sigma,\gamma_{\mathrm{out}} \varphi \rangle_{H^{-1/2}({\Gamma}_\mathrm{out};\R^d),H^{1/2}({\Gamma}_\mathrm{out};\R^d)}\\
    &\qquad\qquad-\langle \gamma_{\mathrm{in}}\varphi,w_1\rangle_{H_0^{1/2}(\Gamma_\mathrm{in};\R^d), H^{-1/2}(\Gamma_\mathrm{in};\R^d)}-\langle \gamma_{\mathrm{out}}\varphi,w_2\rangle_{H_0^{1/2}(\Gamma_\mathrm{out};\R^d)}\\
    &\qquad {=}-\langle \gamma_{\mathrm{in}} \varphi,\calR_{\mathrm{in}}^{-1}\gamma_{\mathrm{n},\mathrm{in}}\sigma+\calR_{\mathrm{in}}^{-1}w_1 \rangle_{H^{1/2}({\Gamma}_\mathrm{in};\R^d)}
    -\langle \gamma_{\mathrm{out}} \varphi,w_2+\calR_{\mathrm{out}}^{-1}u_\sigma \rangle_{H^{1/2}({\Gamma}_\mathrm{out};\R^d)}.
\end{align*}
Since the latter holds for all $\varphi\in H^1_{\sigma,{\mathrm{w}}}(\Omega;\R^d)$, we obtain that
        \begin{align*}
            w_1&=-\gamma_{\mathrm{in}}\sigma-\calR_{\mathrm{in}}u_v,\\
            w_2&=-\gamma_{\mathrm{out}}v-\calR_{\mathrm{out}}^{-1}u_\sigma.
          \end{align*}
These two equations together with \eqref{eq:sysnodeeq1} yield
$(v,u_v,u_\sigma)\in\dom M$ satisfies
\eqref{eq:maxdiss}, which shows the claim.

		\noindent {\em Step~3:} We show that $M$ is closed.\\
		By steps~1\&2 $M$ is dissipative, and $M-\calR$ is onto. Since this implies that $M-\calR$ has closed range, \cite[Prop.~3.14]{EngeNage2000} implies that $M$ is closed.
		\item[\eqref{def:dissnode2}] We show that $P_\calX M= F\&G$ is closed.  As $M$ is closed by \eqref{def:dissnode1}, and in view of \cite[Lem.~2.3]{PhilReis23}, we may equivalently show 
      there exists some $C>0$, such that for all $\spvekk{v}{u_v}{u_\sigma} \in \dom M$
		\begin{align*}
		&\phantom{=}\left\|K\&L \spvekk{v}{u_v}{u_\sigma}\right\|_{H^{-1/2}(\Gamma_\mathrm{in};\R^d)\times H_0^{1/2}(\Gamma_\mathrm{out};\R^d)} \\&\leq C\left( \left\|\spvekk{v}{u_v}{u_\sigma}\right\|_{L^2(\Omega;\R^d) \times H^{1/2}_0(\Gamma_\mathrm{in};\R^d)\times H^{-1/2}(\Gamma_\mathrm{out};\R^d)} + \left\|F\&G\spvekk{v}{u_v}{u_\sigma}\right\|_{L^2(\Omega;\R^d)}\right).
		\end{align*}
	    To this end consider a sequence
		\begin{align}\label{eq:proof_rhsconverges}
        \begin{split}
        \left(\spvekk{v_n}{u_{v,n}}{u_{\sigma,n}}\right)&\quad\text{ bounded in $L^2(\Omega;\R^d)\times H_0^{1/2}(\Gamma_\mathrm{in};\R^d)\times H^{-1/2}(\Gamma_\mathrm{out};\R^d)$},\\
		\left(F\&G\spvekk{v_n}{u_{v,n}}{u_{\sigma,n}}\right)&\quad\text{ bounded in $L^2(\Omega;\R^d)$.}
        \end{split}
		\end{align}
		{\em Step~1:} We show that $(v_n)$ is a~bounded sequence in $H^1_{\sigma,\mathrm{w}}(\Omega;\R^d)$. Since the trace operator $\gamma_{\mathrm{in}}:H^1_{\sigma,\mathrm{w}}(\Omega;\R^d)\to H^{1/2}_0(\Gamma_\mathrm{in};\R^d)$ is bounded and surjective, it possesses a~bounded right inverse \[\gamma_\mathrm{in}^- \in L(H^{1/2}_0(\Gamma_\mathrm{in};\R^d),H^1_{\sigma,\mathrm{w}}(\Omega;\R^d))\] with \[\gamma_{\mathrm{in}}\gamma_\mathrm{in}^-=\Id_{H^{1/2}_0(\Gamma_\mathrm{in};\R^d)}.\]
		This together with the bounded emdedding $H^1_{\sigma,\mathrm{w}}(\Omega;\R^d) \subset L^2(\Omega;\R^d)$ implies that the sequence $(v_n-\gamma^-_{\mathrm{in}} u_{v,n})$ is bounded in $L^2(\Omega;\R^d)$. As $(v_n,u_{v,n},u_{\sigma,n})\in \dom M$, we have that $\gamma_{\mathrm{in}}(v_n-\gamma^-_{\mathrm{in}}u_{v,n})=0$ for all $n\in\mathbb{N}$. Further, let $(P_n)$ be such that \[\mu\nabla v_n -P_n\Id_d\in H(\div,\Omega;\R^{d\times d}).\] Then
		\begin{align*}
		&\phantom{=}
		\langle v_n-\gamma^-_\mathrm{in}u_{v,n},F\&G \spvekk{v_n}{u_{v,n}}{u_{\sigma,n}}\rangle_{L^2(\Omega;\R^d)}\\
        &= \langle v_n-\gamma^-_\mathrm{in}u_{v,n}, \div(\mu \nabla v_n-P_n\Id_d)\rangle_{L^2(\Omega;\R^d)} + \langle v_n-\gamma^-_\mathrm{in}u_{v,n}, \rho(b\cdot \nabla)v\rangle_{L^2(\Omega;\R^d)}\\
		&\underset{\&\eqref{eq:wdiv0}}{\overset{\eqref{eq:normaltracehat}}{=}}-\mu \langle \nabla (v_n-\gamma^-_\mathrm{in}u_{v,n}),\nabla v_n\rangle_{L^2(\Omega;\R^{d \times d})}+\langle v_n-\gamma^-_{\mathrm{D},\Gamma_\mathrm{D}}u_{\mathrm{D},n}, \rho(b\cdot\nabla) v_n\rangle_{L^{2}(\Omega;\R^d)}  \\
        & \qquad \qquad +\langle\gamma_\mathrm{out} (v_n-\gamma^-_{\mathrm{in}}u_{v,n}), \gamma_{\mathrm{n},\mathrm{out}}(\mu \nabla v_n - P_n\Id_d)\rangle_{H^{1/2}_0(\Gamma_\mathrm{out};\R^d) \times H^{-1/2}(\Gamma_\mathrm{out};\R^d)}\\
        &\qquad \qquad +\underbrace{\langle\gamma_\mathrm{in} (v_n-\gamma^-_{\mathrm{in}}u_{v,n}), \gamma_{\mathrm{n},\mathrm{in}}(\mu \nabla v_n - P_n\Id_d)\rangle_{H^{1/2}_0(\Gamma_\mathrm{in};\R^d) \times H^{-1/2}(\Gamma_\mathrm{in};\R^d)}}_{=0} \\
		&=-\mu \langle \nabla v_n,\nabla v_n\rangle_{L^2(\Omega;\R^{d \times d})}+\mu\langle \nabla\gamma^-_\mathrm{in}u_{v,n},\nabla v_n\rangle_{L^2(\Omega;\R^{d \times d})} + \langle v_n,\rho(b\cdot \nabla) v_n\rangle_{L^2(\Omega;\R^d)}\\
        &\qquad \qquad -\langle \gamma^-_\mathrm{in}u_{v,n},\rho(b\cdot\nabla) v_n\rangle_{L^2(\Omega;\R^d)}+ \langle\gamma_\mathrm{out} v_n,u_{\sigma,n}\rangle_{H^{1/2}_0(\Gamma_\mathrm{out};\R^d),H^{-1/2}(\Gamma_\mathrm{out};\R^d)} \\
        &\qquad \qquad -\langle\gamma_\mathrm{out} \gamma^-_\mathrm{in}u_{v,n}, u_{\sigma,n}\rangle_{H^{1/2}_0(\Gamma_\mathrm{out};\R^d),H^{-1/2}(\Gamma_\mathrm{out};\R^d)}\\
		&\leq-\mu \| \nabla v_n\|_{L^2(\Omega;\R^{d \times d})}^2+\mu \|\nabla\gamma^-_\mathrm{in} u_{v,n}\|_{L^2(\Omega;\R^{d \times d})} \|\nabla v_n\|_{L^2(\Omega;\R^{d \times d})} \\
        &\qquad \qquad+ \rho\|v_n\|_{L^2(\Omega;\R^d)}\,\|b\|_{L^\infty(\Omega;\R^d)}\|\nabla v_n\|_{L^2(\Omega;\R^{d \times d})}\\
        &\qquad \qquad+ \rho\| \gamma^-_\mathrm{in}u_{v,n}\|_{L^2(\Omega;\R^d)}\,\|b\|_{L^\infty(\Omega;\R^d)}\|\nabla v_n\|_{L^2(\Omega;\R^{d \times d})}\\
        &\qquad \qquad +\|\gamma_\mathrm{out} v_n\|_{H_0^{1/2}(\Gamma_\mathrm{out};\R^d)}\|u_{\sigma,n}\|_{H^{-1/2}(\Gamma_\mathrm{out};\R^d)}\\&\qquad \qquad + \|\gamma_\mathrm{out} \gamma^-_\mathrm{in}u_{v,n}\|_{H_0^{1/2}(\Gamma_\mathrm{out};\R^d)}\|u_{\sigma,n}\|_{H^{-1/2}(\Gamma_\mathrm{out};\R^d)}
		\end{align*}
		Rearranging the terms in the above inequality, invoking~\eqref{eq:proof_rhsconverges} and continuity of the trace operators, there are constants  $c_1,c_2>0$, such that, for all $n\in\mathbb{N}$,
		\[\| \nabla v_n\|_{L^2(\Omega;\R^{d \times d})}^2\leq c_1+c_2\, \| v_n\|_{H^1(\Omega;\R^d)}.\]
		Thus, as $(v_n)$ is bounded in $L^2(\Omega;\R^d)$, there is $c_3\geq 0$ such that, for all $n\in \mathbb{N}$,
         \[\| v_n\|_{H^1(\Omega;\R^d)}^2\leq c_3+c_2\, \| v_n\|_{H^1(\Omega;\R^d)}.\]
        Applying Young's inequality, this shows that $(v_n)$ is a~bounded sequence in $H_{\sigma,\mathrm{w}}^1(\Omega;\R^d)$. \\
		{\em Step~2:} We show that 
  	    \[
        \begin{split}
        \left(K\&L\spvekk{v_n}{u_{v,n}}{u_{\sigma,n}}\right)&\quad\text{is bounded in $H_0^{1/2}(\Gamma_\mathrm{in};\R^d)\times H^{-1/2}(\Gamma_\mathrm{out};\R^d)$}.
        \end{split}
		\]
By using the Banach-Steinhaus theorem \cite[Thm.~7.3]{Alt16}, it suffices to show that, for all $y = (y_1,y_2)\in H_0^{1/2}(\Gamma_\mathrm{in};\R^d)\times H^{-1/2}(\Gamma_\mathrm{out};\R^d)$, the scalar sequence
		\begin{align*}
		&\left(\left\langle y,K\&L\spvekk{v_n}{u_{v,n}}{u_{\sigma,n}} \right\rangle_{H_0^{1/2}(\Gamma_\mathrm{in};\R^d)\times H^{-1/2}(\Gamma_\mathrm{out};\R^d),H^{-1/2}(\Gamma_\mathrm{out};\R^d)\times H_0^{1/2}(\Gamma_\mathrm{in};\R^d) }\right)\\
        & = \left( \langle y_1,\gamma_{\mathrm{n},\mathrm{in}} (\mu \nabla v_n-P_n\Id_d)\rangle_{H^{1/2}_0(\Gamma_\mathrm{in};\R^d),H^{-1/2}(\Gamma_\mathrm{in};\R^d)} + \langle y_2,\gamma_\mathrm{out} v_n\rangle_{H^{-1/2}(\Gamma_\mathrm{out};\R^d),H^{1/2}_0(\Gamma_\mathrm{out};\R^d)}\right)
		\end{align*}
		is bounded, where, again, $(P_n)$ is such that $\mu\nabla v_n -P_n\Id_d\in H(\div,\Omega;\R^{d\times d})$. As $(v_n)$ is bounded in $H_{\sigma,\mathrm{w}}^1(\Omega;\R^d)$ by Step 1, the sequence \[\left(\langle y_2,\gamma_\mathrm{out} v_n\rangle_{H^{-1/2}(\Gamma_\mathrm{out};\R^d),H^{1/2}_0(\Gamma_\mathrm{out};\R^d)}\right)\] is bounded by continuity of the trace operator $\gamma_\mathrm{out}:H_\mathrm{w}^1(\Omega;\R^d) \to H^{1/2}_0(\Gamma_\mathrm{out};\R^d)$ as defined in~\eqref{eq:rest_dir_trace}.\\		
		As $y_1 \in H^{1/2}_0(\Gamma_\mathrm{in};\R^d)$, we define $w= \gamma_\mathrm{in}^-y_1 \in H^1_{\sigma,\mathrm{w}}(\Omega;\R^d)$, where $\gamma_\mathrm{in}^-$ is the right-inverse of the trace operator as defined in the beginning of Step 1. 
        Then
		\begin{align*}\label{eq:intermediate}
		\begin{split}
		&\phantom{=}\langle w,F\&G \spvekk{v_n}{u_{v,n}}{u_{\sigma,n}}\rangle_{L^2(\Omega;\R^d)} \\&=  -\mu \langle \nabla w,\nabla v_n\rangle_{L^2(\Omega;\R^{d \times d})}+\left\langle y_1,\gamma_{\mathrm{n},\mathrm{in}} (\mu \nabla v_n-P_n\Id_d) \right\rangle_{H^{1/2}_0(\Gamma_\mathrm{in};\R^d),H^{-1/2}(\Gamma_\mathrm{in};\R^d)}\\
        &\qquad + \langle \gamma_\mathrm{out}w,\underbrace{\gamma_{\mathrm{n},\mathrm{out}} (\mu \nabla v_n -P_n\Id_d)}_{=u_{\sigma,n}}\rangle_{H^{1/2}_0(\Gamma_\mathrm{in};\R^d),H^{-1/2}(\Gamma_\mathrm{in};\R^d)} + \langle w, \rho(b\cdot\nabla) v_n\rangle_{L^{2}(\Omega;\R^d)}.
		\end{split}
		\end{align*}
        Now, the scalar sequence $\left\langle y_1,\gamma_{\mathrm{n},\mathrm{in}} (\mu \nabla v_n-P_n\Id_d) \right\rangle_{H^{1/2}_0(\Gamma_\mathrm{in};\R^d),H^{-1/2}(\Gamma_\mathrm{in};\R^d)}$ is bounded, as all other terms in the above identity are bounded.
		\item[\eqref{def:dissnode3}] We show that for all $(u_v,u_\sigma)\in H^{1/2}_0(\Gamma_\mathrm{in};\R^d) \times H^{-1/2}(\Gamma_\mathrm{out};\R^d)$, there exists some $v\in H^1_{\sigma,\mathrm{w}}(\Omega;\R^d)$ with $\spvekk{v}{u_v}{u_\sigma}\in \dom M$. Let $(u_v,u_\sigma)\in H^{1/2}_0(\Gamma_\mathrm{in};\R^d) \times H^{-1/2}(\Gamma_\mathrm{out};\R^d)$. Set $v_{\mathrm{in}}:= \gamma_\mathrm{in}^- u_v \in H^1_{\sigma,\mathrm{w}}(\Omega;\R^d)$, where $\gamma_\mathrm{in}^-$ is the right-inverse defined in Step 1 of (b). 
        Due to the Lax-Milgram lemma, there exists some $v_0\in H^1_{\sigma,\mathrm{w}}(\Omega;\R^d)$ with $\gamma_\mathrm{in} v = 0$ such that, for all $\varphi\in H^1_{\sigma,\mathrm{w}}(\Omega;\R^d)$ with $\gamma_\mathrm{in} \varphi = 0$,
		\begin{align*}
		&\mu\langle  \nabla\varphi,\nabla v_0\rangle_{L^2(\Omega;\R^{d \times d})}+\langle \varphi, v_0\rangle_{L^2(\Omega;\R^d)}-\langle \varphi, \rho(b\cdot\nabla) v_0\rangle_{L^2(\Omega;\R^d)} \\=
		&-\mu\langle \nabla\varphi,\nabla v_\mathrm{in}\rangle_{L^2(\Omega;\R^{d \times d})}-\langle \varphi, v_\mathrm{in}\rangle_{L^2(\Omega;\R^d)}+\langle \varphi, \rho(b\cdot \nabla) v_\mathrm{in}\rangle_{L^2(\Omega;\R^d)}\\
		& -\mu\langle {\gamma_\mathrm{out}\varphi, u_\sigma}\rangle_{H_0^{1/2}(\Gamma_\mathrm{out};\R^d),H^{-1/2}(\Gamma_\mathrm{out};\R^d)}.
		\end{align*}
		By linearity, $v:=v_0+v_\mathrm{in}$ fulfills $\gamma_\mathrm{in} v=u_v$,  $\gamma_{\mathrm{n},\mathrm{out}} (\mu \nabla v-P\Id_d)=u_\sigma$, and, for all $\varphi\in H^1_0(\Omega;\R^d)$ with $\div \varphi = 0$,
		\[
		-\mu \langle \nabla\varphi,\nabla v\rangle_{L^2(\Omega;\R^{d \times d})}+\langle \varphi, \rho(b\cdot\nabla) v\rangle_{L^2(\Omega)}=\langle \varphi, v\rangle_{L^2(\Omega;\R^d)}.
		\]
        Again, by Lemma~\ref{lem:pressure_recovery}, there exists some $P\in L^2(\Omega)$ with $\div \mu \nabla v + \rho(b\cdot\nabla)v - v = \nabla P$ and thus \[\div(\mu\nabla v - P\Id_d) = -\rho(b\cdot\nabla)v + v \in L^2(\Omega;\R^3).\] This gives $\mu \nabla v - P\Id_d\in H(\div,\Omega;\R^{d\times d})$. Hence, we conclude that $\spvekk{v}{u_v}{u_\sigma}\in\dom F\&G$.
		\item[\eqref{def:dissnode4}] We show that $F-\Id$ is surjective, i.e., in particular has dense range. Let $z\in L^2(\Omega;\R^3)$. Again using the Lax-Milgram lemma, there exists some $v\in H^1_{\sigma,\mathrm{w}}(\Omega;\R^d)$ with $\gamma_{\mathrm{in}}v=0$, such that for all $\varphi\in  H^1_{\sigma,\mathrm{w}}(\Omega;\R^d)$ with $\gamma_{\mathrm{in}}\varphi=0$,
		\begin{align}\label{eq:laxmil:ende}
		\mu\langle \nabla\varphi,\nabla v\rangle_{L^2(\Omega;\R^{d \times d})}+\langle \varphi, v\rangle_{L^2(\Omega;\R^d)} -\langle \varphi, \rho(b\cdot \nabla)  v\rangle_{L^2(\Omega;\R^d)}=
		-\langle \varphi,z\rangle_{L^2(\Omega;\R^d)}.
		\end{align}
        Again invoking Lemma~\ref{lem:pressure_recovery}, there exists some $P\in L^2(\Omega)$ with \[\div \mu \nabla v -v+\rho(b\cdot\nabla)v -z = \nabla P,\] which gives rise to 
        \begin{align}\label{eq:divdrauf}
            \div(\mu \nabla v - P\Id_d) =-\rho(b\cdot\nabla)v + z\in L^2(\Omega;\R^d).
        \end{align} Moreover, by $\gamma_\mathrm{in}v = 0$, the Green's identity~\eqref{eq:normaltracehat} together with \eqref{eq:laxmil:ende} and \eqref{eq:divdrauf}, we obtain that \[\gamma_\mathrm{n,out} (\mu \nabla v  - P\Id_d)=0.\] This gives $v\in \dom F$ with $(F-\Id)v = z$.
	\end{itemize}
\end{proof}

\noindent By having shown that $M=\sbvek{F\&G}{K\&L}$ as in \eqref{eq:dissnode_oseen} defines a dissipation node, and by further incorporating the (extremely simple)  Hamiltonian as in \eqref{eq:hamiltonian}, we have brought the Oseen equations 
\begin{align*}
    \dot{p}(t)&=\mu\Delta v(t)+\rho (b\cdot\nabla)v+ \nabla P(t),\!\!\!\!\!\!\!\!\!\!\!\!\!\!\!\!\!\!\!\!\!\!\!\!\!\!\!\!\!\!\!\!\!\!\!\!\!\!\!\!\!\!\!\!\!\!\!\!\!\!\!\!\!\!\!\!\!\!\!\!\!\!\!\!\!\!\!\!\!\!\!\!\!\!\!\!\!\!\!\!\!\!\!\!\!\!\!\!\!\!\!\!\!\!\!\!\!\!\!\!\!\!\!\!\!\!\!\!\!\!\!\!\!\!\!\!\!\!\!\!\!\!\!\!&\\
    \div v(t)&=0,\\
     v\vert_{{\Gamma}_\mathrm{w}}&=0,\\
     v\vert_{{\Gamma}_\mathrm{in}}&=u_v(t),\!\!\!\!\!\!\!\!\!\!\!\!\!\!\!\!\!\!\!\!\!\!\!\!\!\!\!\!\!\!\!\!&\sigma(t)\mathrm{n}\vert_{{\Gamma}_\mathrm{out}}&=u_\sigma(t),\\
     y_\sigma(t)&=\sigma(t)\mathrm{n}\vert_{{\Gamma}_\mathrm{in}},\!\!\!\!\!\!\!\!\!\!\!\!\!\!\!\!\!\!\!\!\!\!\!\!\!\!\!\!\!\!\!\!&
     y_v(t)&=v\vert_{{\Gamma}_\mathrm{out}},\\
     p(0)&=p_0
\end{align*}
with $v(t)=\tfrac1\rho{p(t)}\in L^2(\Omega;\R^d)$ and $\sigma(t)=\mu \left(\nabla v(t) + \nabla v(t)^\top\right) -P(t)\Id_d\in L^2(\Omega;\R^{d\times d})$ into the framework of port-Hamiltonian system nodes. The input consists of velocities at the inflow part $\Gamma_\mathrm{in}$ of the boundary together with the stress tensor in normal direction at the outflow part $\Gamma_\mathrm{out}$ of the boundary. The output is composed of the stress tensor in normal direction at  $\Gamma_\mathrm{in}$ and the velocity trace at $\Gamma_\mathrm{in}$. We have shown that this system is -- under the assumptions on the domain as specified in Section~\ref{sec:domain}, together with $\mu,\rho>0$ and essential boundedness, divergence-freeness and trivial normal boundary trace of $b$  -- port-Hamiltonian in the sense of Definition~\ref{def:pHnode}. This allows to apply the results known for port-Hamiltonian systems of this type, such as, e.g., Proposition~\ref{thm:pHnode}: For instance, we can conclude that the free dynamics of the above system (i.e., $u_v=0$ and $u_\sigma=0$) is described by a~contractive semigroup. We can further conclude from Proposition~\ref{thm:pHnode} that the above system has a~classical solution, if 
\begin{align*}
    u_v&\in W^{2,1}([0,T];H^{1/2}_0(\Gamma_\mathrm{in};\R^d)),\\
u_\sigma&\in W^{2,1}([0,T];H^{-1/2}(\Gamma_\mathrm{out};\R^d)),
\end{align*} and the initial infinitesimal momentum $p_0\in H^1_{\sigma,\mathrm{w}}(\Omega;\R^d)$ has the property 
 that there exists some $P_0\in H^{1}_0(\Omega)$ with $\mu\nabla v-P_0\Id_d\in H(\div,\Omega;\R^{d\times d})$, joint with the compatibility conditions 
\begin{align*}
\tfrac1\rho p_0\vert_{{\Gamma}_\mathrm{in}}&=u_v(0),\\
    \mu \left(\left(\nabla p_0 + \nabla p_0^\top\right) -P_0\Id_d\right)\mathrm{n}\vert_{{\Gamma}_\mathrm{out}}&=u_\sigma(0).
\end{align*}
 Moreover, for the Hamiltonian $\mathcal{H}$ as in \eqref{eq:hamiltonian}, and $0\leq t_0\leq t_1\leq T$,
the weak (and thus also the classical) solutions on the interval $[0,T]$ fulfill the energy balance %
\begin{align*}
 &\phantom{=}\mathcal{H}(p(t_1))-\mathcal{H}(p(t_0))\\ &=  -\mu\int_{t_0}^{t_1}\|\nabla v(t)\|^2_{L^2(\Omega;\R^{d \times d})}\mathrm{d}t  \\&\qquad+  \int_{t_0}^{t_1}\langle y_\sigma,u_{v}\rangle_{H^{-1/2}(\Gamma_\mathrm{in};\R^d),H^{1/2}_0(\Gamma_\mathrm{in};\R^d)}\mathrm{d}t+ \int_{t_0}^{t_1}\langle y_v,u_\sigma\rangle_{H^{-1/2}(\Gamma_\mathrm{out};\R^d),H^{1/2}_0(\Gamma_\mathrm{out};\R^d)}\mathrm{d}t.
\end{align*}

\section{Conclusion}
We have formulated Oseen flows by means of port-Hamiltonian system nodes, providing a functional analytic framework for a system-theoretic and energy-based modeling of boundary controlled linear flow problems. A system node corresponding to the Oseen system has been introduced, and it has been shown that this indeed defines a dissipation node, such that, together with the corresponding kinetic energy Hamiltonian, we obtained an energy balance, linking the change of energy to the enstrophy and the inner product between input and output. Further, we have provided an application to a flow in a tube. 
 
\section*{Acknowledgement}
The authors would like to thank Volker Mehrmann for pointing out the problem of port-Hamiltonian formulation of Oseen flows. 
\bibliographystyle{abbrv}
\bibliography{references.bib}
\end{document}